\documentclass[12pt]{amsart}
\usepackage{graphicx,amssymb,epsfig}
\usepackage{amssymb,amsmath}
\usepackage[all]{xy}
\usepackage{amssymb}
\usepackage{mathrsfs}

\newtheorem{theorem}{Theorem}[section]
\newtheorem{lemma}[theorem]{Lemma}
\newtheorem{proposition}[theorem]{Proposition}
\newtheorem{corollaire}[theorem]{Corollary}

\theoremstyle{definition}
\newtheorem{definition}[theorem]{Definition}

\theoremstyle{remark}
\newtheorem{remark}[theorem]{Remark}

\numberwithin{equation}{section}

\def\R{{\mathbb R}}
\def\Q{{\mathbb Q}}

\def\Z{{\mathbb Z}}

\def\TcT{({\mathbf T} (t))_{t\geq 0}}

\def\TcT0{({\mathbf T}_0 (t))_{t\geq 0}}

\def\L1{L^1 (\R_+ )}

\begin{document}

\begin{center}

{\Large\bf    Stable Homotopy Groups  of Moore Spaces}\\[8mm]

In\`es Saihi \\
Published in \textit{Abualrub, Taher (ed.) et al., Mathematics across contemporary sciences. AUS-ICMS, American University of Sharjah, United Arab Emirates, April 2–5, 2015. Cham: Springer (ISBN 978-3-319-46309-4/hbk; 978-3-319-46310-0/ebook). Springer Proceedings in Mathematics \& Statistics 190, 177-192 (2017)}. 

\end{center}

\bigskip

\begin{center}
\begin{minipage}{24cc}
{\sc Abstract :} 
We determine explicitly the stable homotopy groups of Moore spaces up to the range $7$, using an 
equivalence of categories which allows to consider each Moore space as an exact couple of $\Z$-modules.
\\[5mm]
\end{minipage}
\end{center}


\section{Introduction}

Moore spaces and their stable homotopy groups were widely studied and
a complete reference on this subject is the book of Baues \cite{Baues}.

In this paper, we propose a new approach allowing to see Moore spaces as exact couples of
$\Z$-modules by means of an equivalence of categories. Even though a similar result is proven in
\cite{Baues}, the approach given here is of independent interest, since it is used to determine
explicitly the stable homotopy groups of Moore spaces up to the range $7$.

 Let $G$ be an abelian group and $n$ an integer greater than $1$. A Moore space $M(G\mathbin,n)$
 is a simply connected CW-complex $X$ such that $H_n(X)\simeq G$ and 
 $\widetilde{ H_i}(X)=0$ for $i\not=n$. The homotopy type of $M(G\mathbin, n)$ is uniquely
determined by the pair $(G\mathbin, n)$ (see \cite{Vogel1}). 

Let ${\mathscr M}_n$ be the category whose objects are Moore spaces
$M(A\mathbin, n)$, where $A$ is a $\Z$-module,  and whose morphisms are 
homotopy classes  of pointed maps between such Moore spaces. 
Notice that, unlike the Eilenberg-MacLane, 
the set homotopy classes of pointed maps $[M(A\mathbin, n)\mathbin,M(B\mathbin, n)]$ 
between two Moore spaces is different from ${\rm Hom}(A\mathbin,B)$ (see proposition \ref{prop}).

Let ${\mathcal Mod}$ be the category of $\Z$-modules and
let ${\mathscr D}_e$ be the category of exact couples in ${\mathcal Mod}$ 
$$
\xymatrix{ A \ar[rr]^{2}  && A\ar[ld]^{\alpha} \\ & B\ar[lu]^{\beta}}
$$
such that $\alpha\beta=2$.

There are two exact functors $\Phi_1$ and $\Phi_2$ from ${\mathscr D}_e$ 
to ${\mathcal Mod}$ assigning to a diagram the $\Z$-module $A$ or $B$ respectively.

The aim of  section $2$ is to construct, for $n\geq3$, 
an equivalence of categories $\mathcal E$ between ${\mathscr M}_n$  and ${\mathcal D}_e$.  In
\cite{Baues} and in a different context, Baues gave a similar result using the properties of the
 Whitehead $\Gamma$-functor.

In  section $3$, the stable homotopy groups $\pi^S_{i}(X)$ ($0\leq i\leq7$) of a Moore space $X$ will be 
expressed in terms of $\Phi_j({\mathcal E}(X))$ for $j=1,2$. The same techniques can be used 
to determine $\pi_i^S(M(A\mathbin,n))$ for $i\geq8$, but calculations become complicated.

\section{Equivalence of Categories between Moore Spaces and Diagrams} \label{secmain}

\subsection{Category of Diagrams } In this section, we propose an equivalence of categories 
that allows to consider Moore spaces as diagrams of $\Z$-modules.

 Recall that the suspension functor from ${\mathscr M}_{n}$ to ${\mathscr M}_{n+1}$ is an equivalence of 
 categories for $n\geq3$, so next results are independent of $n$.

Consider two modules $A$ et $B$.
Let $X$ be the Moore space $X=M(A\mathbin,n)$, $Y$ the Moore space $Y=M(B\mathbin,n)$ and
 $[X\mathbin,Y]$  the set of  homotopy classes of pointed maps from $X$ to $Y$; this set is an abelian group 
 (see \cite{Hatcher}). Moreover:

\begin{proposition} \label{prop}(\cite{Baues}, \cite{Hatcher})\quad
There is a natural exact sequence:
\begin{equation}\label{exactseq}
\xymatrix{\relax
   0\ar[r]^-{} & {\rm Ext}(A\mathbin,B/2)\ar[r]^-{} &[X\mathbin,Y]\ar[r]^-{} & {\rm Hom}(A\mathbin,B) \ar[r]^-{} &0.}
 \end{equation}
\end{proposition}
\bigskip

Set $S=M({\Z}\mathbin,n)$ and  $P=M({\Z}/2\mathbin,n)$. Applying the exact sequence \eqref{exactseq}
to  $S$ and $X$, we obtain:
$$
\xymatrix{\relax
   0\ar[r]^-{} & {\rm Ext}(\Z\mathbin,A/2)\ar[r]^-{} &[S\mathbin,X]\ar[r]^-{} & {\rm Hom}(\Z\mathbin,A) \ar[r]^-{} &0}
  $$
hence $[S\mathbin,X]\simeq A$. Similarly, applied to $P$ and $X$,  \eqref{exactseq} becomes:
  $$
\xymatrix{\relax
0\ar[r]^-{} & {\rm Ext}(\Z/2\mathbin,A/2)\ar[r]^-{} &[P\mathbin,X]\ar[r]^-{} & {\rm Hom}(\Z/2\mathbin,A)
\ar[r]^-{} &0.}
$$
Since ${\rm Ext}(\Z/2\mathbin,A/2)$ is naturally isomorphic to $A/2$ (see proposition \ref{isomZ/2mod}),
and if $A'=[P\mathbin,X]$, then:
\begin{equation}\label{Equ2}
\xymatrix{\relax
   0\ar[r]^-{} & A/2\ar[r]^-{} &A'\ar[r]^-{} & A_2 \ar[r]^-{} &0}
 \end{equation}
where $A_2$ denotes order $2$ elements of $A$. In other words, we have the long exact sequence:
\begin{equation}\label{Equ1}
\xymatrix{\relax
   A\ar[r]^-{2} &A\ar[r]^-{\alpha_X} &A' \ar[r]^-{\beta_X} &A \ar[r]^-{2} &A
   } 
   \end{equation}
or equivalently, the exact couple denoted by $D_X$:
 $$\xymatrix{ A \ar[rr]^{2}  && A\ar[ld]^{\alpha_X} \\ & A'\ar[lu]^{\beta_X}}$$
 
 Moreover,
if $f$ is a map betwen two Moore spaces $X=M(A\mathbin,n)$ and 
$Y=M(B\mathbin,n)$, then we can deduce a map
 $\bar f : D_X\longrightarrow D_Y$ as follows:
$\bar f =(f_1\mathbin,f_2)$ where $f_1:A\simeq[S\mathbin, X]\longrightarrow B\simeq[S\mathbin, Y]$  and 
$f_2:A'=[P\mathbin, X]\longrightarrow B'=[P\mathbin, Y]$ 
are the natural maps induced by $f$. The following diagrams commute:
$$
\xymatrix{\relax
   A\ar[r]^-{2}\ar[d]^{}_{f_1}&A\ar[r]^-{\alpha_X} \ar[d]^{}_{f_1}&A' \ar[r]^-{\beta_X} \ar[d]^{f_2}_{}&A \ar[r]^-{2} \ar[d]^{f_1}_{}&A\ar[d]^{f_1}_{}\\
  B\ar[r]^{2} &B\ar[r]^-{\alpha_Y} &B' \ar[r]^-{\beta_Y} &B \ar[r]^-{2} &B  
   } $$%%%%%%%

\subsubsection{Particular case of $P$}

When $X=P$, we have the next results:

\begin{proposition} 
$ [P\mathbin,P]\simeq{\Z}/4$.
\end{proposition} 
\noindent The proof of this result can be found in \cite{Hatcher}  or \cite{Wu}.

\begin{lemma}\label{CasP}
The composition $\alpha_P\beta_P $ is multiplication by $2$ on  ${\Z}/4$.
\end{lemma}
\noindent 
{\it Proof}\quad When $X=P$, the exact sequence \eqref{Equ1} becomes:
$$
\xymatrix{\relax
   \Z/2\ar[r]^-{2} &\Z/2\ar[r]^-{\alpha_P} &\Z/4 \ar[r]^-{\beta_P} &\Z/2 \ar[r]^-{2} &\Z/2
   } $$%%%%%%% 
i.e: 
   $$
\xymatrix{\relax
0\ar[r]^-{ } &\Z/2\ar[r]^-{\alpha_P} &\Z/4 \ar[r]^-{\beta_P} &\Z/2 \ar[r]^-{} &0
   } $$%%%%%%%  
then $\alpha_P=2$  and $\beta_P$ is the canonical surjection. Hence
$\alpha_P\beta_P\mathbin:{\Z}/4\buildrel{2}\over\longrightarrow {\Z}/4$.

\subsubsection{General Case}

\begin{lemma}
For any Moore space $X=M(A\mathbin, n)$, the composition
$\alpha_X\beta_X\mathbin:A'\buildrel{}\over\longrightarrow A'$ is multiplication by $2$.
\end{lemma}  

\noindent 
{\it Proof}\quad Let $u\in A'=[P\mathbin,X]$ and $f\mathbin:P\longrightarrow X$ a representative 
of  $u$. We have two maps $f_1$ and $f_2$ and the commutative diagram: 
$$
\xymatrix{\relax
   \Z/2\ar[r]^-{2=0}\ar[d]^{}_{f_1}&\Z/2\ar[r]^-{\alpha_P} \ar[d]^{}_{f_1}&\Z/4 \ar[r]^-{\beta_P} \ar[d]^{f_2}_{}&\Z/2 \ar[r]^-{2=0} \ar[d]^{f_1}_{}&\Z/2\ar[d]^{f_1}_{}\\
  A\ar[r]^-{2} &A\ar[r]^-{\alpha_X} &A' \ar[r]^-{\beta_X} &A \ar[r]^-{2} &A  
   } $$%%%%%%%  
If $u_0$  denotes the class of the
identity map in  $[P\mathbin,P]$, then $f_2(u_0)=u$. The result is an immediate consequence of
 lemma \ref{CasP}.

\subsubsection{Category of diagrams}
\begin{definition} 
Let ${\mathscr D}_e$ be the category of exact couples in the category ${\mathcal Mod}$ of $\Z$-modules
 \begin{equation}\label{diag}
\xymatrix{ A \ar[rr]^{2}  && A\ar[ld]^{\alpha} \\ & B\ar[lu]^{\beta}}
\end{equation}
such that $\alpha\beta=2$.

A  morphism $f$ between two objetcs $D$ and $D'$ is a couple $f=(f_1\mathbin,f_2)$
such that the following diagrams commute: 
$$
\xymatrix{\relax
   A\ar[r]^-{2}\ar[d]^{}_{f_1}&A\ar[r]^-{\alpha} \ar[d]^{}_{f_1}&B \ar[r]^-{\beta} \ar[d]^{f_2}_{}&A \ar[r]^-{2} \ar[d]^{f_1}_{}&A\ar[d]^{f_1}_{}\\
  A'\ar[r]^{2} &A'\ar[r]^-{\alpha'} &B' \ar[r]^-{\beta'} &A' \ar[r]^-{2} &A' 
   } $$%%%%%%%   
\end{definition}   

%\begin{remark}
%Cette cat\'egorie ${\mathcal D}$ n'est pas ab\'elienne. Par contre, si l'on consid\`ere ${\mathcal D}_0$ la 
%cat\'egorie dont les objets sont les diagrammes commutatifs 
%$$\xymatrix{ A \ar[rr]^{2}  && A\ar[ld]^{\alpha} \\ & B\ar[lu]^{\beta}}$$
%""Let ${\mathscr D}$ be the category of diagrams in the category ${\mathcal Mod}$
% \begin{equation}\label{diag}
%\xymatrix{ A \ar[rr]^{2}  && A\ar[ld]^{\alpha} \\ & B\ar[lu]^{\beta}}
%\end{equation}
  %%
%such that
%$2\alpha=2\beta=\beta\alpha=0$ and $\alpha\beta=2$.""
  %%
%o\`u $A$ et $B$ sont deux modules et avec les conditions $2\alpha=0$, $2\beta=0$,
%$\beta\alpha=0$ et $\alpha\beta=2$, et dont 
% les morphismes sont identiques aux morphismes de ${\mathcal D}$, alors
%  ${\mathcal D}_0$  est une cat\'egorie ab\'elienne.
%\end{remark}
\medskip
{\it Notations :}\quad For ease, an object of ${\mathscr D}_e$ will be denoted by
$$
\xymatrix{ A \ar@<2pt>[r] ^{\alpha}& B \ar@<2pt>[l]^{\beta}}$$
%%%%%%%%
and a morphism between two objects $D$ and $D'$ will be denoted by
  $$
  \xymatrix{
   A  \ar@<2pt>[r] ^{\alpha_X} \ar[d]^{}_{f_1} & 
  A' \ \ar@<2pt>[l] ^{\beta_X}\ar[d]^{f_2}_{}  \\
    **[l]B  \ar@<2pt>[r] ^{\alpha_Y}& 
    **[r]B'  \ar@<2pt>[l] ^{\beta_Y} }
    $$
%%%%%%

\medskip The previous constructions can be summarized in the following statement:
\begin{proposition}
There is a functor $\mathcal E:{\mathscr M}_n\longrightarrow{\mathscr D}_e$
assigning to each Moore space $X$ the diagram $D_X$, and to each homotopy class $f$ of
pointed maps between two Moore spaces $X$ and $Y$ the map $\bar f : D_X\longrightarrow D_Y$.
\end{proposition}

In the remaining of this section, we will prove that the functor ${\mathcal E}$ is an equivalence of
categories.

\medskip {\it Notations :}\quad Let $\Phi_1$ and $\Phi_2$ denote the two functors from ${\mathscr D}_e$ 
to ${\mathcal M}od$ defined as follows:
if $D$ is an object of${\mathscr D}$ given by:
\begin{equation}\label{diagexact}
\xymatrix{ A \ar@<2pt>[r] ^{\alpha}& B \ar@<2pt>[l]^{\beta}}
\end{equation}
%%%%%%%%
then $\Phi_1(D)=A$ and $\Phi_2(D)=B$.

Notice that there is a natural transformation between functors $\Phi_1$ and $\Phi_2$ obtained  
by associating to a diagram $D$ given by \eqref{diagexact}, the morphism $\alpha$. By associating 
to the diagram $D$ the morphism $\beta$, we get a natural transformation from $\Phi_2$ to $\Phi_1$.

\subsection{Equivalence of categories between ${\mathscr M}_n$ and ${\mathscr D}_e$}
\subsubsection{Some Algebraic Results} This section is devoted to prove some general
 algebraic results needed  to obtain the equivalence of categories announced above.
\begin{proposition}\label{isomZ/2mod}
\quad
For every ${\Z}/2$-modules $A$ and $B$, there is an isomorphism $\lambda_{(A\mathbin,B)}$, natural in
$A$ and in $B$:
$$\xymatrix{\relax
\lambda_{(A\mathbin,B)}: {\rm Ext}(A\mathbin,B)\ar[r]^-{\sim}& {\rm Hom}(A\mathbin,B).}$$
\end{proposition}

\noindent 
{\it Proof :}
An element $e$ of ${\rm Ext}(A\mathbin,B)$ is represented by an extension: 
$$
\xymatrix{\relax
   0\ar[r]^-{} & B\ar[r]^-{f} &E\ar[r]^-{g} & A \ar[r]^-{} &0.}
$$
%%%%%%%
Each element $x\in A$ is of order $2$ and $g$ is surjective, so there is $y\in E$ such that $g(y)=x$
and $2y\in\ker g={\rm Im}f$. Since $f$ is injective, there exists a unique $z\in B$ such that $f(z)=2y$.
The map assigning to $x$ the element $z$ is well defined; then we obtain a morphism:
$$
\xymatrix{\relax
\lambda_{(A\mathbin,B)}: {\rm Ext}(A\mathbin,B)\ar[r]^-{}& {\rm Hom}(A\mathbin,B).}$$
Since $A$ is free, there is a natural isomorphism
${\rm Ext}(A\mathbin,B)\longrightarrow{\rm Hom}(A\mathbin,{\rm Ext}({}\Z/2\mathbin,B))$
 obtained by restriction. 
 (Each $a\in A$ defines a map ${\Z}/2\longrightarrow A$ which induices an extension of ${\Z}/2$ 
 by $B$ using a pull-back.)
But ${\rm Ext}({}\Z/2\mathbin,B)$ is naturaly isomorphic to $B$, 
so we get an isomorphism from ${\rm Ext}(A\mathbin,B)$ to ${\rm Hom}(A\mathbin,B)$, 
which is $\lambda_{(A\mathbin,B)}$.

\medskip
\begin{remark}
If $A$ and $B$ are two ${\Z}$-modules, we construct similarly a  morphisme natural in $A$ and $B$, 
denoted also by  $\lambda_{(A\mathbin,B)}$  : 
$$\lambda_{(A\mathbin,B)}\mathbin:{\rm Ext}(A\mathbin,B)\longrightarrow{\rm Hom }(A_2\mathbin,B/2)$$
and obtained by the composition:
$$
\xymatrix{\relax
\lambda_{(A\mathbin,B)}: {\rm Ext}(A\mathbin,B)\ar[r]^-{} &{\rm Ext}(A_2\mathbin,B/2)
\ar[r]^-{\lambda_{(A_2\mathbin,B/2)}} &{\rm Hom }(A_2\mathbin,B/2),
} $$%%%%%%%  
where the first  morphism is induced by restriction to order  $2$ elements in $A$ 
and the projection of $B$ on $B/2$.
\end{remark}

\medskip
\begin{corollaire}
If $A$ is a ${\Z}/2$-module and $B$ a ${\Z}$-module, then
${\rm Ext}(A\mathbin,B)\simeq{\rm Hom }(A\mathbin,B/2)$.
\end{corollaire}

\noindent 
{\it Proof :}  %Comme $A$ est un ${\Z}/2$-module, il est libre et on peut l'\'ecrire comme
%$A=\oplus{\Z}/2$. Par suite, ${\rm Ext}(A\mathbin,B)
%\simeq{\rm Hom}(A\mathbin, {\rm Ext}({\Z}/2\mathbin,B))$.
%Comme ${\rm Ext}({\Z}/2\mathbin,B)\simeq B/2$, on obtient le r\'esultat.
The morphism $\lambda_{(A\mathbin,B)}$ is the composition %(??? pourquoi?)
$$
\xymatrix{\relax
  \lambda_{(A\mathbin,B)}: {\rm Ext}(A\mathbin,B)\ar[r]^-{ pr} &{\rm Ext}(A\mathbin,B/2)
  \ar[r]^-{\lambda_{(A\mathbin,B/2)}} &{\rm Hom }(A\mathbin,B/2)
   } $$%%%%%%%  
where ${pr}$ is the morphism induced by the projection of $B$ on $B/2$.
By \ref{isomZ/2mod},  $\lambda_{(A\mathbin,B/2)}$ is an isomorphism; it suffices to show that
${ pr}\mathbin:{\rm Ext}(A\mathbin,B)\longrightarrow
{\rm Ext}(A\mathbin,B/2)$ is bijective. 
%(???)
But $A$ is a ${\Z}/2$-module, so $A$ is free and then can be written $A=\oplus{\Z}/2$ . Since 
${\rm Ext}(\oplus{\Z}/2\mathbin,B)=\prod{\rm Ext}({\Z}/2\mathbin,B)$ we can show the result for $A=\Z/2$.
Using the resolution 
$$
\xymatrix{\relax
   0\ar[r]^-{} &{\Z}\ar[r]^-{2} &{\Z} \ar[r]^-{} &{\Z}/2 \ar[r]^-{} &0,
   } $$%%%%%%%  
we get the diagram:
$$
\xymatrix{\relax
  {\rm Hom}({\Z}/2\mathbin,B)\ar[r]^-{}\ar[d]^{}_{\simeq}& B\ar[r]^-{2}\ar[d]^{}_{pr}&B\ar[r]^-{} \ar[d]^{}_{pr}& {\rm Ext}({\Z}/2\mathbin,B) \ar[r]^-{} \ar[d]^{\simeq}_{}&0 \\
B_2\ar[r]^-{}&B/2\ar[r]^-{2=0}& B/2\ar[r]^-{} & {\rm Ext}({\Z}/2\mathbin,B/2) \ar[r]^-{} &0
   } $$%%%%%%%  
   \medskip
\begin{corollaire}
If $A$ is a $\Z$-module and $B$ a ${\Z}/2$-module, then
${\rm Ext}(A\mathbin,B)\simeq{\rm Hom }(A_2\mathbin,B)$.
\end{corollaire}   

\noindent 
{\it Proof :} Since the morphism $\lambda_{(A\mathbin,B)}$ is the composition %(??? pourquoi?)
$$
\xymatrix{\relax
  \lambda_{(A\mathbin,B)}: {\rm Ext}(A\mathbin,B)\ar[r]^-{\rm R} &{\rm Ext}(A_2\mathbin,B)
  \ar[r]^-{\lambda_{(A_2\mathbin,B)}} &{\rm Hom }(A_2\mathbin,B),
   } $$%%%%%%% 
where $R$ is the morphism induced by the restriction to $A_2$, and ${\lambda_{(A_2\mathbin,B)}}$ is
an isomorphism, we have just to show that $R$ is bijective.

But $B$ is free, so $B=\oplus{\Z}/2$; consider the injective module $I=\oplus (\Q/\Z)$; then
we have the exact sequence:
$$
\xymatrix{\relax
   0\ar[r]^-{} & B\ar[r]^-{} &I\ar[r]^-{2} & I \ar[r]^-{} &0.}
$$
Applying the functor ${\rm Hom}(A\mathbin,\cdot)$,  we obtain the following diagram:
$$
\xymatrix{\relax
{\rm Hom}(A\mathbin,I)\ar[r]^-{2}\ar[d]^{}_{R}& {\rm Hom}(A\mathbin,I)\ar[r]^-{} \ar[d]^{}_{R}&
{\rm Ext}(A\mathbin,B) \ar[r]^-{} \ar[d]^{R}_{}&0 \\
{\rm Hom}(A_2\mathbin,I)\ar[r]^-{2=0}& {\rm Hom}(A_2\mathbin,I)\ar[r]^-{} &
{\rm Ext}(A_2\mathbin,B) \ar[r]^-{} &0
} $$%%%%%%%  
where $R$ denotes the morphism induced by the restriction to $A_2$.

On the other hand, we have the exact sequence: 
 $$
\xymatrix{\relax
   0\ar[r]^-{} &A_2\ar[r]^-{} &A\ar[r]^-{2} &A }
$$
%%%%%%%
Applying the functor ${\rm Hom}(\cdot\mathbin,I)$, we get an isomorphism between 
${\rm Hom}(A_2\mathbin,I)$ and ${\rm Hom}(A\mathbin,I)/2$, so
$R\mathbin:{\rm Ext}(A\mathbin,B)\longrightarrow{\rm Ext}(A_2\mathbin,B)$ is bijective.

\medskip
\subsubsection{Equivalence of categories}

\begin{theorem}\label{equivalence}
The functor $\mathcal E$  is an equivalence of categories between
${\mathscr M}_n$ and ${\mathscr D}_e$.
\end{theorem}

To prove this theorem, we need 
the next two lemmas.
\begin{lemma}\label{objets}
For each diagram $D$ in ${\mathscr D}_e$, there exists a Moore space $X$ in
${\mathscr M}_n$ such ${\mathcal E}(X)=D$.
\end{lemma} 

\noindent 
{\it Proof:}
Let $D$ be an object of ${\mathscr D}_e$ given by:
 $$
 \xymatrix{ A \ar@<2pt>[r] ^{\alpha}& B \ar@<2pt>[l]^{\beta}}$$
%%%%%%%%%%%%%%%%%%%%%
Set $X$ the Moore space
$X=M(A\mathbin,n)$. The diagram associated to $X$ is given by:
$$
\xymatrix{ A \ar@<2pt>[r] ^{\alpha_X}& A' \ar@<2pt>[l]^{\beta_X}}$$
%%%%%%%%%%%%%%%%
Then, we have the following diagram:
$$
\xymatrix{\relax
   0\ar[r]^-{} & A/2\ar[r]^-{\alpha} \ar[d]^{}_{Id}&B \ar[r]^-{\beta} &A_2 \ar[r]^-{}\ar[d]^{}^{Id} &0\\
 0\ar[r]^-{} & A/2\ar[r]^-{\alpha_X} &A' \ar[r]^-{\beta_X} &A_2 \ar[r]^-{} &0   }
$$
Each horizontal exact sequence defines an element in ${\rm Ext}(A_2\mathbin,A/2)\simeq 
{\rm Hom}(A_2\mathbin,A/2)$. Since $\beta\alpha =2$ on $B$ and
$\beta_X\alpha_X=2$ on $A'$, the two extensions give the same element in
${\rm Hom}(A_2\mathbin,A/2)$ and then the two extensions are idsomorphic.

\medskip

\begin{lemma}\label{morphismes}
If $X$ and $Y$ are two Moore spaces, then
$$
\xymatrix{\relax
 [X\mathbin,Y]\ar[r]^-{\simeq} &{\rm Hom }(D_X\mathbin,D_Y) =
 {\rm Hom }({\mathcal E}(X)\mathbin,{\mathcal E}(Y))   } $$%%%%%%% 
\end{lemma}

\noindent 
{\it Proof}\quad Let $X=M(A\mathbin, n)$ and $Y=M(B\mathbin, n)$, then there is an exact sequence:
$$
\xymatrix{\relax
   0\ar[r]^-{} & {\rm Ext}(A\mathbin,B)\ar[r]^-{} &[X\mathbin,Y]\ar[r]^-{} & {\rm Hom}(A\mathbin,B) \ar[r]^-{} &0}
    $$
But $ {\rm Ext}(A\mathbin,B/2)\simeq {\rm Ext}(A_2\mathbin,B/2)\simeq  {\rm Hom}(A_2\mathbin,B/2)$, 
so we obtain the exact sequence:
$$
\xymatrix{\relax
 0\ar[r]^-{} & {\rm Hom}(A_2\mathbin,B/2)\ar[r]^-{} &[X\mathbin,Y]\ar[r]^-{} & 
{\rm Hom}(A\mathbin,B) \ar[r]^-{} &0.}
$$
%%%
On the other side, the forgetful morphism
${Fr}\mathbin :{\rm Hom }(D_X\mathbin,D_Y)\longrightarrow{\rm Hom}(A\mathbin,B)$ is surjective.
Recall that an element  $g\in{\rm Hom }(D_X\mathbin,D_Y)$ is given by two maps $g_1$ and $g_2$
such that:
$$
\xymatrix{\relax
A\ar[r]^-{2}\ar[d]^{}_{g_1}&A\ar[r]^-{\alpha_X} \ar[d]^{}_{g_1}&A' \ar[r]^-{\beta_X} \ar[d]^{g_2}_{}&A 
\ar[r]^-{2} \ar[d]^{g_1}_{}&A\ar[d]^{g_1}_{}\\
B\ar[r]^{2} &B\ar[r]^-{\alpha_Y} &B' \ar[r]^-{\beta_Y} &B \ar[r]^-{2} &B  
   } $$%%%
so an element $g\in{\rm Hom }(D_X\mathbin,D_Y)$ is in the kernel of the forgetful morphism if $g_1=0$
and then we obtain a morphism $A_2\longrightarrow B/2$.
Hence, we get the following commutative diagram:
$$
\xymatrix{\relax
  0\ar[r]^-{} & {\rm Hom}(A_2\mathbin,B/2)\ar[r]^-{} \ar[d]^{}_{f_{X\mathbin,Y}}&[X\mathbin,Y]\ar[r]^-{} \ar[d]^{}_{}
  & {\rm Hom}(A\mathbin,B) \ar[r]^-{} \ar[d]^{Id}_{}&0\\
0\ar[r]^-{} & {\rm Hom}(A_2\mathbin,B/2)\ar[r]^-{} &{\rm Hom}(D_X\mathbin,D_Y)\ar[r]^-{Fr} & 
{\rm Hom}(A\mathbin,B) \ar[r]^-{} &0   }
$$
%%%%
To prove the isomorphism between $[X\mathbin,Y]$ and ${\rm Hom }(D_X\mathbin,D_Y)$, 
it suffices to  verify that $f_{X\mathbin,Y}$ is the identity map. 
Notice that $f_{X\mathbin,Y}$ is a bifunctor, covariant in $B$ and contravariant in $A$.

When $X=P$ and $Y=S$, the diagram becomes:
$$
\xymatrix{\relax
  0\ar[r]^-{} & \Z/2\ar[r]^-{} \ar[d]^{}_{f_{P\mathbin,S}}&[P\mathbin,S]\simeq\Z/2\ar[r]^-{} \ar[d]^{}_{}& 0 \ar[r]^-{} &0\\
0\ar[r]^-{} & \Z/2\ar[r]^-{} &{\rm Hom}(D_P\mathbin,D_S)\simeq\Z/2\ar[r]^-{Fr} &0 \ar[r]^-{} &0   }
$$
%%%%
so $f_{P\mathbin,S}$ is necessarily  the identity map. 
When $X=P$ and $Y=M(B\mathbin,n)$: an element $y\in B$ defines a morphism $\Z\longrightarrow B$ 
that can be realized by a map between Moore spaces  $S\longrightarrow Y$.
Let ${ \bar y}\in {\rm Hom}({\Z}/2\mathbin, B/2)$ denote the quotient of the map defined by $y$.
By assigning ${ \bar y}$ to the generator of $\Z/2$, we get the commutative diagram:
\begin{equation*}%\label{diagfonctoriel}
\xymatrix{\relax
 {\Z}/2\ar[r]^-{} \ar[d]^{}_{f_{P\mathbin,S}=Id}&{\rm Hom}({\Z}/2\mathbin,B/2)\simeq B/2
 \ar[d]^{f_{P\mathbin,Y}}_{}\\
 {\Z}/2\ar[r]^-{} &{\rm Hom}({\Z}/2\mathbin,B/2)\simeq B/2   }
\end{equation*}
%%%%
Since ${\rm Hom}({\Z}/2\mathbin, B/2)$ is naturally  isomorphic to $ B/2$, 
even in this case $f_{P\mathbin,Y}=Id$.

Given $x\in A_2$, it defines a map ${\Z}/2\longrightarrow A_2\subset A$ which can be realized by
a map of Moore spaces $P\longrightarrow X$. This map allows to have the following commutative
diagram, using the functoriality of $f_{X\mathbin,Y}$:
$$
\xymatrix{\relax
 {\rm Hom}(A_2\mathbin,B/2)\ar[r]^-{} \ar[d]^{}_{f_{X\mathbin,Y}}&{\rm Hom}({\Z}/2\mathbin,B/2)\simeq 
 B/2\ar[d]^{f_{P\mathbin,Y}=Id}_{}\\
  {\rm Hom}(A_2\mathbin,B/2)\ar[r]^-{} &{\rm Hom}({\Z}/2\mathbin,B/2)\simeq B/2   }
$$
the horizontal maps assign to a morphism $\varphi\mathbin:A_2\longrightarrow B/2$ its evaluation
$\varphi(x)\in B/2$.
To conclude that $f_{X\mathbin,Y}$ is the identity map on $ {\rm Hom}(A_2\mathbin,B/2)$, it suffices
to notice that the module $A_2$ is $\Z/2$-free, and if  $\{u_i\}_{ i\in I}$ is a basis of $A_2$ then
${\rm Hom}(A_2\mathbin,B/2) \simeq\prod{\rm Hom}({\Z}/2\mathbin,B/2)\simeq\prod B/2$. Using the
evaluation on each generator $u_i$, we deduce the desired result.

\begin{remark}
With lemmas \ref{objets} et \ref{morphismes}, we get the proof of theorem \ref{equivalence}.
\end{remark}
\section{Stable Homotopy Groups of Moore Spaces}  

\bigskip

Let  $X=M(A\mathbin,n)$ and consider the Atiyah-Hirzebruch spectral sequence in homology
with coefficients in the satble homotopy groups:
$$  H_p(X\mathbin; \pi_q^S)\Rightarrow\pi_{p+q}^S(X).$$
This spectral sequence contains just two non trivial columns  and induces the following exact sequence:
\begin{equation}\label{SSAH}
\xymatrix{\relax
   0\ar[r]^-{} & A\otimes\pi_{q}^{S}\ar[r]^-{\nu^X} &\pi_{n+q}^{S}(X) \ar[r]^-{\mu^X} & {\rm Tor}(A\mathbin,\pi_{q-1}^{S}) \ar[r]^-{} &0}
   \end{equation}
    %%%%%%%
Moreover, this exact sequence is natural in $X$.

Notice that, if $\underline{X}$ denotes the spectrum associated to the Moore  space $X$,
then  $\pi_{n+i}^S(X)=\pi_{i}^S(\underline{X})$. 
In the following, the spectrum associated to a space $X$ will also be denoted  by $X$.

Recall the first stable homotopy groups (see  \cite{Koch}):
$$\pi_0^S=\Z\mathbin, \pi_1^S={\Z}/2 \mathbin, \pi_2^S={\Z}/2
\mathbin, \pi_3^S={\Z}/24 \mathbin, \pi_4^S=\pi_5^S=0 \mathbin, \pi_6^S={\Z}/2
\mathbin, \pi_7^S={\Z}/240,$$
so the exact sequence \eqref{SSAH} allows to obtain, for any Moore space $X$:
$$\pi_0^S(X)=A,\qquad \pi_1^S(X)\simeq A\otimes{\Z}/2=A/2,\qquad 
\pi_4^S(X)\simeq {\rm Tor}(A\mathbin,{\Z}/24)=A_{24}\mathbin,$$
$$\pi_5^S(X)=0
,\qquad\pi_6^S(X)\simeq A\otimes{\Z}/2=A/2 $$
but we can't determine  explicitly $\pi_2^S(X)$ , $\pi_3^S(X)$ and $\pi_7^S(X)$.

To compute $\pi_i^S(X)$, for $i=2\mathbin,3\mathbin,7$, we need the following lemma:

\medskip
\begin{lemma}
$\pi_2^S(P)={\Z}/4$, $\pi_3^S(P)={\Z}/2\oplus{\Z}/2$, $\pi_7^S(P)={\Z}/2\oplus{\Z}/2$.
\end{lemma}

\noindent 
{\it Proof}\quad These groups are given in \cite{Wu}, but we propose an easer proof of these results
using the  arguments of section $2$.

For $q=2\mathbin,3\mathbin, 7$, the exact sequence \eqref{SSAH} becomes:
$$
\xymatrix{\relax
0\ar[r]^-{} & \Z/2\ar[r]^-{} &\pi_{q}^{S}(P) \ar[r]^-{} &\Z/2 \ar[r]^-{} &0}
$$
then $\pi_q^S(P)\simeq {\Z}/2\oplus{\Z}/2$ or $\pi_q^S(P)\simeq{\Z}/4$.

There is a cofibration sequence:
\begin{equation}\label{cofibration}
\xymatrix{\relax
   \ar[r]^-{\theta} & P\ar[r]^-{\delta} &S \ar[r]^-{2} &S \ar[r]^-{\theta} &P\ar[r]^-{\delta} &}
   \end{equation}
where $\theta$ is of degree $0$ and $\delta$ of degree $-1$. If $\lambda$ denotes
 the composition of $\delta$ by the Hopf map from $S$ to $S$, then we get the Moore spectra 
 diagram
\begin{equation}\label{DiagMoore}
 \xymatrix{ S \ar[rr]^{2}  && S\ar[ld]^{\theta} \\ & P\ar[lu]^{\lambda}}
 \end{equation}
verifying
 $2\theta=0$, $2\lambda=0$, $\lambda\theta=0$ et $\theta\lambda=2$. 

Applying the functor $\pi_2^S$ to \eqref{DiagMoore}, we obtain the following diagram :
 $$
 \xymatrix{ \Z/2 \ar[rr]^{2=0}  && \Z/2\ar[ld]^{\theta_*} \\ & \pi_2^S(P)\ar[lu]^{\lambda_*}}$$
where $2\lambda_*=0$, $2\theta_*=0$ and $\theta_*\lambda_*=2$.  
This diagram is not necessarily exact, but, the exact sequence of stable homotopy groups applied to the
cofibration \eqref{cofibration}  gives:
$$\ker(\theta_*:{\Z}/2\longrightarrow\pi_2^S(P))={\rm Im(2:{\Z}/2\longrightarrow{\Z}/2}).$$
% car  si l'on $X\longrightarrow Y\longrightarrow Z$ ...... alors.......[H].....????.

Then it suffices to find an element $u\in\pi_2^S(P)$ such that $\lambda_*(u)=1$. For this purpose, we can
choose $n=2$ so $S=S^2$ and $P=P_2= \Sigma{\R}P_2$. 
We have the cofibration:
$$S^2\longrightarrow P_2\longrightarrow S^3.$$ 
Applying the stable homotopy functor, we get:
$$
\xymatrix{\relax
  \pi_4^S(S^2)\ar[r]^-{}\ar@{=}[d]& \pi_4^S(P_2)\ar[r]^-{}\ar@{=}[d]& \pi_4^S(S^3) \ar[r]^-{} \ar@{=}[d]& 
  \pi_3^S(S^2) \ar[r]^-{} \ar@{=}[d]& \pi_3^S(P_2)\ar@{=}[d]\ar[r]^-{}& \pi_3^S(S^3)\ar@{=}[d]\\
  \pi_2^S \ar[r]^-{}\ar@{=}[d]& \pi_2^S(P)\ar[r]^-{}\ar@{=}[d]& \pi_1^S \ar[r]^-{} \ar@{=}[d]& \pi_1^S \ar[r]^-{} 
  \ar@{=}[d]& \pi_1^S(P)\ar@{=}[d]\ar[r]^-{}& \pi_0^S\ar@{=}[d]\\
  \Z/2\ar[r]{} &\pi_2^S(P)\ar[r]^-{}  \ar@{->>}[rd]&\Z/2 \ar[r]^-{0}\ar[d]_{\simeq}^{\rm Hopf} &\Z/2
   \ar[r]^-{Id} &\Z/2\ar[r]^{0} &\Z\\
  &&\Z/2&&& &
   } $$%%%%%%%   
Then $\pi_2^S(P)$ is surjected on ${\Z}/2=\pi_1^S$ which is sent by the Hopf map on ${\Z}/2=\pi_2^S$
by assigning to the generator $\eta$ of $\pi_1^S={\Z}/2$ the  generator $\eta^2$ of $\pi_2^S={\Z}/2$.

Now, applying the functor $\pi_3^S$ to \eqref{DiagMoore}, we get:
  $$
  \xymatrix{ \Z/24 \ar[rr]^{2}  && \Z/24\ar[ld]^{\theta_*} \\ & \pi_3^S(P)\ar[lu]^{\lambda_*}}$$
with $2\theta_*=0$, $2\lambda_*=0$ and $\theta_*\lambda_*=2$.  This diagram is not necessarily exact, 
but
  $$\ker(\theta_*:{\Z}/24\longrightarrow\pi_3^S(P))={\rm Im(2:{\Z}/24\longrightarrow{\Z}/24}).$$

Let $x\in\pi_3^S(P)$, then $2\lambda_*(x)=0$. There exists $u\in{\Z}/24$ such that $\lambda_*(x)=12u$.
So  $2x=\theta_*(\lambda_*(x))=\theta_*(12u)=0$ since $2\theta_*=0$. This implies that 
elements of $\pi_3^S(P)$ vanish when multiplied by $2$ and then 
$\pi_3^S(P)\simeq {\Z}/2\oplus{\Z}/2$.

The same argument shows that $\pi_7^S(P)={\Z}/2\oplus{\Z}/2$.

\bigskip

Consider a Moore space $X=M(A\mathbin,n)$. The next theorems compute $\pi_i^S(X)$, for $i=2,3,7$, 
in terms of the modules $\Phi_1(D_X)$ and $\Phi_2(D_X)$.

\begin{theorem}\label{pi2stable}
For each generator $\gamma\in\pi_2^S(P)$, there is a natural isomorphism
$\pi_2^S(X)\simeq \Phi_2(D_X)=[P\mathbin,X]$.
\end{theorem}
 
\noindent 
{\it Proof }\quad Consider the exact sequence \eqref{Equ2} and the exact sequence \eqref{SSAH} for $q=2$ :
$$
\xymatrix{\relax
   0\ar[r]^-{} & A/2\ar[r]^-{\alpha} &A'\ar[r]^-{\beta} & A_2 \ar[r]^-{} &0}
$$
%%%%%%%    

$$
\xymatrix{\relax
   0\ar[r]^-{} & A/2\ar[r]^-{\nu^X} &\pi_2^S(X)\ar[r]^-{\mu^X} &A_2 \ar[r]^-{} &0}
$$
%%%%%%%  
We construct a map $A'\longrightarrow \pi_2^S(X)$ as follows:
choose $\gamma$ a generator of $\pi_2^S(P)\simeq{\Z}/4$.
Let $u\in A'$ and consider $f$  representing the class $u\in A'=[P\mathbin,X]$.Then  $f$ induces a
map $f_*\mathbin: \pi_2^S(P)\longrightarrow  \pi_2^S(X)$ and we define
 $\varphi_{\gamma}(u)=f_*(\gamma)$. the map  $\varphi_{\gamma}$ relies the two exact sequences:
 $$
\xymatrix{\relax
   0\ar[r]^-{} & A/2\ar[r]^-{\alpha} &A'\ar[r]^-{\beta}\ar[d]^{\varphi_{\gamma}}_{} & A_2 \ar[r]^-{} &0\\
   0\ar[r]^-{} & A/2\ar[r]^-{\nu^X} &\pi_2^S(X)\ar[r]^-{\mu^X} & A_2 \ar[r]^-{} &0 }
$$
%%%%%%% 
Now, we may prove that the composite map
$$
\xymatrix{\relax
   A'\ar[r]^-{\varphi_{\gamma}} &\pi_2^S(X)\ar[r]^-{\mu^X} & A_2 }
$$
is $\beta\mathbin: A'\longrightarrow A_2$.
Using the functoriality, it suffices to prove the result when $X=P$.
In this case, he diagram becomes: 
$$
\xymatrix{\relax
   0\ar[r]^-{} & {\Z}/2\ar[r]^-{\alpha} &{\Z}/4\ar[r]^-{\beta}\ar[d]^{\varphi_{\gamma}}_{} & { \Z}/2 
   \ar[r]^-{} \ar[d]^{Id}&0\\
   0\ar[r]^-{} &  {\Z}/2\ar[r]^-{\nu^P} &\pi_2^S(P)\simeq {\Z}/4\ar[r]^-{\mu^P} &  {\Z}/2 \ar[r]^-{} &0 }
$$
%%%%%%%   
where we see clearly that $\mu^P\circ\varphi_{\gamma}=\beta$.

By functoriality, for each Moore space $X=M(A\mathbin,n)$ we get the following commutative diagram:
$$
\xymatrix{\relax
   0\ar[r]^-{} & A/2\ar[r]^-{\alpha} \ar[d]^{}_{h}&A'\ar[r]^-{\beta}\ar[d]^{\varphi_{\gamma}}_{} & A_2 \ar[r]^-{} \ar[d]^{Id}&0\\
   0\ar[r]^-{} & A/2\ar[r]^-{\nu^X} &\pi_2^S(X)\ar[r]^-{\mu^X} &A_2 \ar[r]^-{} &0 }
$$
%%%%%%% 
here $h$ is the natural map making the diagram commute. Notice that $h$ is functorial in $X$. 
Then, to determine  $h\mathbin:A/2\longrightarrow A/2$, it suffices to study the case
$X=S$. In that case, the diagram becomes: 
 $$
\xymatrix{\relax
   0\ar[r]^-{} & {\Z}/2\ar[r]^-{\alpha} \ar[d]^{}_{h}&{\Z}/2\ar[r]^-{\beta}\ar[d]^{\varphi_{\gamma}}_{} & 0 \ar[r]^-{} \ar[d]^{Id}&0\\
   0\ar[r]^-{} & {\Z}/2\ar[r]^-{\nu^S} &\pi_2^S={\Z}/2\ar[r]^-{\mu^S} &0 \ar[r]^-{} &0 }
$$
%%%%%%%
and then $h$ is necessarily the identity map.

Let  $X=M(A\mathbin,n)$ be a Moore space. Each element $x \in A$ defines a maps
$f \mathbin: \Z \longrightarrow A$ given by $f(1)=x$. This map is realized  by a map between Moore spaces 
$f\mathbin:S\longrightarrow X$ and induces, by naturality of $h$, the following commutative diagram:   
 $$
\xymatrix{\relax
   {\Z}/2\ar[r]^-{\bar{f}} \ar[d]^{}_{h=Id}&A/2\ar[d]^{h}_{} \\
   {\Z}/2\ar[r]^-{\bar{f}} &A/2 }
$$
%%%%%%%  
so $h\mathbin:A/2\longrightarrow A/2$ is still the identity map.
\medskip
\begin{remark}
The isomorphism $\pi_2^S(X)\simeq A'$ depends on the choice of the generator
$\gamma\in \pi_2^S(P)={\Z}/4$. Choosing the generator $-\gamma$ multiplies the
isomorphisme by $-1$. 
\end{remark}

\begin{theorem}
For each $\gamma\in\pi_3^S(P)$ such that $\mu^P(\gamma)=1\in{\Z}/2$, there is a natural isomorphism
$\pi_3^S(X)\simeq A'\oplus_{A/2}A/24$ obtained by the pushout
$$
\xymatrix{\relax
  A/2\ar[r]^-{\alpha} \ar[d]^{}_{\times12}& A'\ar[d]^{}_{} \\
  A/24\ar[r]^-{} &\pi_3^S(X) }
$$
%%%%%%% 
where  $A=\Phi_1(D_X)$ and $A'=\Phi_2(D_X)$.
\end{theorem}
\noindent 
{\it Proof}\quad  When $q=3$, the exact sequence \eqref{SSAH} becomes:
$$
\xymatrix{\relax
   0\ar[r]^-{} & A/24\ar[r]^-{\nu^X} &\pi_3^S(X)\ar[r]^-{\mu^X} &A_2 \ar[r]^-{} &0}
$$
%%%%%%%  
For $X=P$ we get 
$$
\xymatrix{\relax
   0\ar[r]^-{} & \Z/2\ar[r]^-{\nu^P} &\pi_3^S(P)\simeq\Z/2\oplus\Z/2\ar[r]^-{\mu^P} &\Z/2 \ar[r]^-{} &0}
$$
%%%%%%% 
Choose $\gamma\in\pi_3^S(P)$ such that $\mu^P(\gamma)$ is the generator of ${\Z}/2$. 
We construct a map
$\varphi_{\gamma}\mathbin: A'\longrightarrow\pi_3^S(X)$ as follows: 

Let $u\in A'=[P\mathbin,X]$ and let $f\mathbin:P\longrightarrow X$ representing the class $u$. Then
$\varphi_{\gamma}(u)=f_*(\gamma)$.

As in  the proof of theorem \ref{pi2stable} we show that the composition of  
$\mu^X:\pi_3^S(X)\longrightarrow A_2$ by $\varphi_{\gamma}$ is $\beta\mathbin: A'\longrightarrow A_2$.

We obtain the following commutative diagram:
 $$
\xymatrix{\relax
   0\ar[r]^-{} & A/2\ar[r]^-{\alpha} \ar[d]^{}_{h}&A'\ar[r]^-{\beta}\ar[d]^{\varphi_{\gamma}}_{} & A_2 
   \ar[r]^-{} \ar[d]^{Id}&0\\
   0\ar[r]^-{} & A/24\ar[r]^-{\nu^X} &\pi_3^S(X)\ar[r]^-{\mu^X} &A_2 \ar[r]^-{} &0 }
$$
%%%%%%% 
Since $h$ is natural, we need just to determine it for $X=S$.
In that case, the map $h\mathbin:{\Z}/2\longrightarrow {\Z}/24$ assigns to the generator of ${\Z}/2$ 
an element of ${\Z}/24$ vanishing when multiplied by $2$, 
that means $0$ or $12$. Then $h=0$ or $h=\times12$.  To prove that $h=\times12$, we 
consider the cofibration 
$$
\xymatrix{\relax
   S \ar[r]^-{2} &S \ar[r]^-{} &P}
    $$
which induces the long exact sequence:
 \begin{equation}\label{SEHS}
\xymatrix{\relax
   \cdots\ar[r]^-{} & \pi_n^S\ar[r]^-{2} &\pi_n^S \ar[r]^-{} &\pi_n^S(P) \ar[r]^-{} &\pi_{n-1}^S\ar[r]^-{2} &\pi_{n-1}^S\ar[r]&\cdots}
 \end{equation}
  For $n=3$,  we have:
$$
\xymatrix{\relax
\pi_3^S={\Z}/24\ar[r]^-{2} &\pi_3^S={\Z}/24 \ar[r]^-{} &\pi_3^S(P) \ar[r]^-{} &\pi_{2}^S={\Z}/2
\ar[r]^-{2=0} &\pi_{2}^S={\Z}/2}
$$
This proves that $\pi_3^S(P)\longrightarrow\pi_{2}^S$ is surjective, so every map $S\longrightarrow S$
 of degree $2$ can be lifted to a map $S\longrightarrow P$ of degree $3$.

If $\gamma\in\pi_3^S(P)$ is represented by  a map, denoted also $\gamma\mathbin: S^5\longrightarrow P_2$, 
and since $\mu^P(\gamma)=1\in{\Z}/2= {\rm Tor}({\Z}/2\mathbin,\pi_{2}^{S}) $, then the map
$\pi_3^S(P)\longrightarrow\pi_{2}^S$ takes $\gamma$ to the generator $1_{{\Z}/2}\in\pi_2^S={\Z}/2$.

Let
$$
\xymatrix{\relax
  u\mathbin: P_2 \ar[r]^-{\delta_2} &S^3 \ar[r]^-{Hopf} &S^2}
 $$
 be a representative of the nonzero element of $[P\mathbin,S]={\Z}/2$
and $a\mathbin:S^5\longrightarrow S^3$ a representative of the generator of $\pi_2^S$. Then
 $\varphi_{\gamma}([u])=u_*(\gamma)=\eta\times(\delta_2)_*(\gamma)=\eta\times[a]$ where
$\eta$ denotes the multiplication by the class of the Hopf map. 
But the multiplication by the Hopf map class takes the generator of $\pi_2^S$ to product by $12$ of
the generator of $\pi_3^S$  (see \cite{Koch}).
This allows to deduce that
 $\varphi_{\gamma}([u])=12\in{\Z}/24$  and that $h$ is multiplication by $12$.

\begin{remark}
 The isomorphism $\pi_3^S(X)\simeq A'\oplus_{A/2}A/24$ depends on the choice of
$\gamma\in \pi_3^S(P)\simeq{\Z}/2\oplus{\Z}/2$ verifying $\mu^P(\gamma)=1$. 
There are two possible choices.

If we choose $\gamma'$ such that $\mu^P(\gamma')=1$, then $\nu^P(1_{{\Z}/2})=\gamma-\gamma'$.
We can show that
$$\varphi_{\gamma'}=\varphi_{\gamma}+\tilde{\lambda}\circ\beta$$
where $\tilde{\lambda}\mathbin:A_2\longrightarrow\pi_3^S(X)$ is defined as follows:
if $a\in A_2$, we can represent it by a map $a\mathbin:S\longrightarrow X$ such that $2a=0$. 
This map induces $a_*\mathbin:\pi_3^S\longrightarrow\pi_3^S(X)$ 
taking all generators of $\pi_3^S$ to the same element $a_*(1_{{\Z}/24})\in\pi_3^S(X)$ since $2a_*=0$. 
Then we define $\tilde\lambda$ by $\tilde{\lambda}(a)=a_*(1_{{\Z}/24})$.
\end{remark}
\medskip

\begin{theorem}
For each $\gamma\in\pi_7^S(P)$ such that $\mu^P(\gamma)=1\in{\Z}/2$, there is a natural isomorphism 
 $\pi_7^S(X)\simeq A/240\oplus A_2 $, where $A=\Phi_1(D_X)$.
\end{theorem}

\noindent 
{\it Proof}\quad Using the same construction of the case of $\pi_3^S(X)$, we get the following
commutative diagram:
 $$
\xymatrix{\relax
   0\ar[r]^-{} & A/2\ar[r]^-{\alpha} \ar[d]^{}_{h}&A'\ar[r]^-{\beta}
   \ar[d]^{\varphi_{\gamma}}_{} & A_2 \ar[r]^-{} \ar[d]^{Id}&0\\
   0\ar[r]^-{} & A/240\ar[r]^-{\nu^X} &\pi_7^S(X)\ar[r]^-{\mu^X} &A_2 \ar[r]^-{} &0 }
$$
%%%%%%% 
To determine $h$, it suffices to consider the case of $X=S$, since it is natural on $X$. In that case, 
$h\mathbin:{\Z}/2\longrightarrow {\Z}/240$, so $h=0$ or $h=\times120$.

For $n=7$, the long exact sequence \eqref{SEHS} becomes:
$$
\xymatrix{\relax
\pi_7^S={\Z}/240\ar[r]^-{2} &\pi_7^S={\Z}/240 \ar[r]^-{} &\pi_7^S(P) \ar[r]^-{} &\pi_{6}^S={\Z}/2\ar[r]^-{2=0} 
&\pi_{6}^S={\Z}/2}
$$
showing that $\pi_7^S(P)\longrightarrow\pi_{6}^S$ is surjective. We use the same techniques of the
previous theorem proof, and the fact that the product by the Hopf class on $\pi_6^S$ is zero
(see \cite{Koch}), we deduce that $h=0$

\begin{remark}
Using the new universal coefficient exact sequence of  \cite{saihi}, we can represent the functor
$\pi_i^S$ on ${\mathscr M}_n$ as a tensor product by particular objects of an abelian category
${\mathscr D}$ containing ${\mathscr D}_e$.
\end{remark}

\bigskip

\bibliographystyle{amsplain}

\begin{thebibliography}{10}

%\bibitem{[B]  } L. Breen, R. Mikhailov \emph{Derived Functors Of Non-additive Functors And Homotopy Theory} Arxiv preprint arXiv : 0910. 2817 2009-arxiv.org

\bibitem{Baues} H-J. Baues, { Homotopy Type And Homology}, Oxford Mathematical Monographs, Clarendon Press, 1996.


\bibitem{Hatcher} A.Hatcher, { Algebraic Topology}, Cambridge University Press, Cambridge 2000.

\bibitem{Koch} S. O. Kochman, {Stable homotopy groups of spheres}, Lecture Notes in Math. no. 1423, (1990)

\bibitem{saihi} I. Saihi,  {\it Homologies g\'en\'eralis\'ees \`a coefficients}, C. R. Acad. Sci. Paris, Ser. I 353 
(2015) 397-401

\bibitem{Switzer}  R. M. Switzer, {Algebraic Topology-Homology and Homotopy}, Springer-Verlag, Berlin 1975.

\bibitem{Vogel1}  P. Vogel, {\it A solution of the Steenrod problem for $G$-Moore spaces}, K-theory 1 (1987), 
no. 4, 325-335.

\bibitem{Wu}  J. Wu, {\it Homotopy Theory of the Suspensions of the Projective Plane} (2003)- books.google.com


\end{thebibliography}

\providecommand{\bysame}{\leavevmode\hbox to3em{\hrulefill}\thinspace}

\vfill
{ines.saihi@fst.utm.tn}\\
University of Tunis El Manar, Faculty of Sciences of Tunis

\end{document}